%% file: main.tex
\documentclass[11pt,a4paper]{article}

\usepackage[T1]{fontenc}
\usepackage[utf8]{inputenc}
\DeclareUnicodeCharacter{2009}{\,}
\usepackage{amsmath}
\usepackage{amssymb}
\usepackage{amsthm}
\usepackage[a4paper,margin=1in]{geometry}
\usepackage{hyperref}
\hypersetup{hidelinks}

\newtheorem{theorem}{Theorem}
\newtheorem{lemma}{Lemma}

\title{Optimisation of Robin Laplacian eigenvalue with indefinite weight in spherical shell}
\author{%
Baruch Schneider$^{1,*}$ \quad Diana Schneiderov\'a$^{1}$ \quad Yifan Zhang$^{1}$\\[0.75em]
\small $^{1}$Department of Mathematics, University of Ostrava, Ostrava, Czech Republic\\
\small $^{*}$Corresponding author: \href{mailto:baruch.schneider@osu.cz}{\texttt{baruch.schneider@osu.cz}}
}
\date{}

\begin{document}

\maketitle

\begin{abstract}
This paper is concerned with an optimisation problem of Robin Laplacian eigenvalue with respect to an indefinite weight, which is formulated as a shape optimisation problem thanks to the known bang-bang distribution of the optimal weight function. The minimisation of the principal eigenvalue of the problem in a spherical shell of an arbitrary dimension is fully solved.
\end{abstract}


\input{sections/introduction}

\input{sections/result}
\input{sections/proof}
\input{sections/acknowledgment}

\bibliographystyle{plainurl}
\bibliography{bibliography/references}

\end{document}

%% file: sections/introduction.tex
\section{Introduction}\label{sec:introduction}

The study of eigenvalue problems in a spherical shell is rich in both mathematical structure and applications.
In quantum mechanics, the eigenvalues of the Laplacian relate to the energy levels of particles confined in potential wells. Robin boundary conditions can model systems with specific physical constraints, such as quantum vacuum fluctuations (\cite{MunozCastaneda2015}) or entanglements (\cite{Robin2021}). 
In acoustics, the eigenvalues of the Laplacian can describe the resonant frequencies of spherical shells (\cite{Fritze2005}), which are relevant in designing musical instruments or architectural acoustics.
In fluid dynamics, the study of eigenvalues can help analyze flow stability (\cite{Clarte2021}) in spherical geometries. 
In electromagnetics, eigenvalue problems can describe the modes of electromagnetic waves in spherical cavities (\cite{Esposito1999}), which is important in antenna design and wireless communication.

The eigenvalue problem considered in this paper is given as
\begin{equation}\label{eq:pb}
	\left\{
		\begin{aligned}
			\Delta \phi + \lambda m\phi  &=0 &\ \text{in} &\ \Omega\\
			\partial_{\mathbf{n}} \phi + \beta\phi &=0 &\ \text{on} &\ \partial\Omega
		\end{aligned}
	\right.
\end{equation}
with Robin boundary condition ($\beta\ge 0$), 
where $\Omega$ is a bounded domain in $\mathbb R^n$ with Lipschitz boundary, $\mathbf{n}$ is the unit outer normal on $\partial\Omega$, and the weight function $m$ is a bounded measurable function that changes sign in $\Omega$, i.e. $\Omega_m^+:=\{x\in\Omega;\ m(x)>0\}$ has measure between $0$ and $|\Omega|$, and that
$$-1\le m(x)\le \kappa\ \text{a.e.}\ x\in \Omega$$
for a constant $\kappa>0$. 

According to \cite{Lou2006}, \cite{Hintermueller2011} and \cite{Lamboley2016}, $\lambda$ is called a \emph{principal eigenvalue} of \eqref{eq:pb}, if the associated eigenfunction $\phi$ is positive, where $\phi\in H^1(\Omega)$. 

The existence and uniqueness of a positive principal eigenvalue is shown in \cite{Afrouzi1999} and \cite{Bocher1914}. It is also well-known (cf. \cite{Brown1980}, \cite{Senn1982}, \cite{Lou2006} and \cite{Lamboley2016}) that the eigenvalue can be expressed with a Rayleigh quotient. To sum up, we have:

\begin{theorem}\label{thm:principal}
If $\beta\ge 0$, the problem \eqref{eq:pb} admits a unique positive principal eigenvalue $\lambda$, depending on $m$. Moreover, 
$$\lambda(m) = \inf_{\phi\in \mathcal{S}(m)} \frac{\int_\Omega |\nabla \phi|^2 + \beta\int_{\partial\Omega}\phi^2}{\int_\Omega m\phi^2},$$
where $\mathcal{S}(m):=\left\{\phi\in H^1(\Omega); \int_\Omega m\phi^2>0\right\}$.
\end{theorem}

We are interested in the minimisation of the positive principal eigenvalue $\lambda(m)$:
\begin{equation}\label{eq:pb-v1}
	\begin{aligned}
		\inf_{m\in \mathcal{M}_{m_0,\kappa}} &\lambda(m),\\
		\mathcal{M}_{m_0,\kappa} &:= \left\{m\in L^\infty(\Omega);\ -1\le m\le \kappa, \left|\Omega_m^+\right|>0, \int_\Omega m\le -m_0|\Omega|\right\}
	\end{aligned}
\end{equation}
Here $m_0$ is a real constant such that $m_0\in(-\kappa, 1)$ if $\beta>0$ and $m_0\in(0,1)$ if $\beta=0$. 

The following result is proven in \cite{Lou2006} for Neumann boundary condition ($\beta=0$), but the extension to Robin case is straightforward:
\begin{theorem}\label{thm:bangbang}
	There exists some weight function $m^*\in 	\mathcal{M}_{m_0,\kappa}$ such that the infimum of \eqref{eq:pb-v1} is attained at $\lambda(m^*)$. Moreover, there exists a measurable subset $E^*$ of $\Omega$ such that 
	\begin{equation}
		m^* = \kappa \chi_{E^*} - \chi_{\Omega\setminus E^*} \text{ a.e. } x\in \Omega,
	\end{equation}
	where $\chi_U$ denotes the characteristic function of a set $U$. In addition, the volume constraint is active at $m^*$:
	$$\int_\Omega m^* = -m_0|\Omega|.$$
\end{theorem}

In other words, the optimal weight function $m^*$ satisfies a bang-bang distribution in $\Omega$ with respect to some optimal set $E^*$. The problem can therefore be formulated as a shape optimisation one:
\begin{equation}\label{eq:pb-v2}
	\begin{aligned}
		\inf_{E\in \mathcal{E}_{c,\kappa}}  \lambda(E) & :=\lambda(\kappa \chi_E - \chi_{\Omega\setminus E}),\\
		\mathcal{E}_{c,\kappa} &:= \{E\subseteq \Omega; E\ \text{is measurable}, |E|\in (0,c|\Omega|]\}
	\end{aligned}
\end{equation}
where $c:=\frac{1-m_0}{1+\kappa}$ and $c\in (0,1)$ if $\beta>0$ while $c\in\left(0,\frac{1}{1+\kappa}\right)$ if $\beta=0$. The optimal $E$ is expressed as a function of $\beta$, $\kappa$ and $c$. 

The 1-dimensional result is fully obtained in \cite{Lamboley2016} as follows: 

\begin{theorem}\label{thm:sol1d}
	Let $n=1$ and $\Omega=(0,1)$, and
	\begin{equation}
		\beta^*=\beta^*(c,\kappa):=\begin{cases}
			\frac{2}{c\sqrt{\kappa}} \arctan\left(\frac{1}{\sqrt{\kappa}}\right)\ &\text{if}\ \kappa>1,\\
			\frac{\pi}{2c}\ &\text{if}\ \kappa=1,\\
			\frac{1}{c\sqrt{\kappa}}\left(\arctan\left(\frac{2\sqrt{\kappa}}{\kappa-1} + \pi\right)\right)\ &\text{if}\ \kappa<1.
		\end{cases}
	\end{equation}
	Then the optimal set $E^*$ for the problem \eqref{eq:pb-v2} is solved as follows:
	\begin{itemize}
		\item if $\beta>\beta^*$, $E^* = \left(\frac{1}{2}-\frac{c}{2},\frac{1}{2}+\frac{c}{2}\right)$ is the unique solution (up to a set of measure zero);
		\item if $\beta=\beta^*$, $E^*$ can be any open interval in $\Omega$ of length $c$;
		\item if $\beta<\beta^*$, $E^*=(0,c)$ or $E^*=(1-c,1)$. 
	\end{itemize}
\end{theorem}

For higher dimensions, a complete result is still unknown. The cylindrical domains are studied in \cite{Kao2008} with Neumann boundary conditions and in \cite{Hintermueller2011} with Robin boundary conditions. Moreover, in \cite{Lamboley2016}, some higher-dimensional case is studied when $\partial\Omega$ is connected, and it is shown that the optimal set is not always a ball. In general, when $n\ge 2$, determination of the optimal shape $E^*$ is widely open. 

In this paper we consider the case when $\Omega$ is a $n$-spherical shell. The case when $\beta=0$ corresponds to Neumann boundary condition and is studied in \cite{Jiang2023}. The idea of changing variables from \cite{Jiang2023} turns out to be useful in the Robin case as well, however, in \cite{Jiang2023}  the accurate scaling of new weight function after the change of variable is ignored, and also the new variable for the case $n\ge 3$ is defined incorrectly as a positive variable.

%% file: sections/result.tex
\section{Results}\label{sec:result}

Denote by $B_\rho(x)$ the open ball centred at $x$ of radius $\rho$, where $x\in \mathbb R^n$ and $\rho\ge 0$. An $n$-spherical shell is given by $A_{r_1,r_2} := B_{r_2}(0)\setminus \overline{B_{r_1}(0)}$ for $r_1,r_2\in\mathbb R$ satisfying $0<r_1<r_2$. Consider the problem \eqref{eq:pb} in a shell $A_{r_1,r_2}$, and we shall prove the following results.

\begin{theorem}\label{thm:shell-2d}
	Let $n=2$, $\Omega = A_{r_1,r_2}$ and denote by $r:=(x_1+x_2)^{\frac{1}{2}}$. Consider the problem \eqref{eq:pb} in $\Omega$, then there exists some constant $m_0'\in (0,1)$ such that the following hold:
	\begin{itemize}
		\item if $\beta >\frac{\beta^*}{r_1(\ln(r_2)-\ln(r_1))} $, the optimal set for $r$ is 
		$$\left(r_1^{\frac{1+c'}{2}} r_2^{\frac{1-c'}{2}}, r_1^{\frac{1-c'}{2}} r_2^{\frac{1+c'}{2}}\right);$$
		\item if $\beta = \frac{\beta^*}{r_1(\ln(r_2)-\ln(r_1))} $, the optimal set for $r$ is
		$$\left(r_0, r_0 r_1^{-c'}r_2^{c'}\right)$$
		for any $r_0\in \left[r_1,r_1^{c'}r_2^{1-c'}\right]$; 
		\item if $\beta <\frac{\beta^*}{r_1(\ln(r_2)-\ln(r_1))} $, the optimal set for $r$ is 
		$$\left(r_1,r_1^{1-c'}r_2^{c'}\right)\ \text{or}\ \left(r_1^{c'}r_2^{1-c'}, r_2\right).$$
	\end{itemize}
	Here $c':=\frac{1-m_0'}{1+\kappa}$, and $\beta^*=\beta^*(c',\kappa)$ is defined as in Theorem \ref{thm:sol1d}. 
\end{theorem}

\begin{theorem}\label{thm:shell-3d}
	Let $n$ be a natural number satisfying $n\ge 3$, $\Omega = A_{r_1,r_2}$ and denote by $r:=(x_1+x_2+\cdots +x_n)^{\frac{1}{2}}$. Consider the problem \eqref{eq:pb} in $\Omega$, then there exists some constant $m_0'\in (0,1)$ such that the following hold:
	\begin{itemize}
		\item if $\beta > \frac{(n-2) r_1^{1-n} \beta^*}{\left(r_1^{2-n} - r_2^{2-n}\right)}$, the optimal set for $r$ is
		$$\left(\left(\frac{(1+c')r_1^{2-n}+(1-c')r_2^{2-n}}{2}\right)^{\frac{1}{2-n}}, \left(\frac{(1-c')r_1^{2-n}+(1+c')r_2^{2-n}}{2}\right)^{\frac{1}{2-n}} \right);$$
		\item if $\beta = \frac{(n-2) r_1^{1-n} \beta^*}{\left(r_1^{2-n} - r_2^{2-n}\right)}$, the optimal set for $r$ is
		$$\left([(2-n)t_0]^{\frac{1}{2-n}}, \left[(2-n)t_0+c'(r_2^{2-n}-r_1^{2-n})\right]^{\frac{1}{2-n}}\right)$$
		for any $t_0\in \left[\frac{r_1^{2-n}}{2-n}, \frac{(1-c')r_2^{2-n}+c'r_1^{2-n}}{2-n}\right]$;
		\item if $\beta < \frac{(n-2) r_1^{1-n} \beta^*}{\left(r_1^{2-n} - r_2^{2-n}\right)}$, the optimal set for $r$ is 
		$$\left(r_1,\left[(1-c')r_1^{2-n}+c'r_2^{2-n}\right]^{\frac{1}{2-n}}\right)\ \text{or}\ \left(\left[(1-c')r_2^{2-n}+c'r_1^{2-n}\right]^{\frac{1}{2-n}},r_2\right).$$
	\end{itemize}
	Here $c':=\frac{1-m_0'}{1+\kappa}$, and $\beta^*=\beta^*(c',\kappa)$ is defined as in Theorem \ref{thm:sol1d}. 
\end{theorem}

%% file: sections/proof.tex
\section{Proof of results}\label{sec:proof}

We first give a straightforward generalisation of Theorem \ref{thm:sol1d} to any interval:

\begin{lemma}\label{lem:1d}
	Let $n=1$ and $\Omega=(a,b)$ for $a,b\in\mathbb R$ satisfying $a<b$. Then the optimal set $E^*$ to the problem \eqref{eq:pb-v2} is given as follows:
	\begin{itemize}
		\item if $\beta>\frac{\beta^*}{b-a}$, 
		$$E^*=\left(\frac{a+b}{2}-\frac{c(b-a)}{2}, \frac{a+b}{2}+\frac{c(b-a)}{2}\right);$$
		\item if $\beta=\frac{\beta^*}{b-a}$, $E^*$ can be any open interval in $\Omega$ of length $c(b-a)$;
		\item if $\beta<\frac{\beta^*}{b-a}$, $E=(a,a+c(b-a))$ or $E=(b-c(b-a),b)$.
	\end{itemize}
	Here the $\beta^*$ is defined the same as in Theorem \ref{thm:sol1d}. 
\end{lemma}
\begin{proof}
	When $x\in(a,b)$, let $y:=\frac{x-a}{b-a}$, then $y\in(0,1)$. The lemma then follows from Theorem \ref{thm:sol1d}. 
\end{proof}

From now on we prove Theorem \ref{thm:shell-2d} and Theorem \ref{thm:shell-3d}. First we switch to the spherical coordinates. Let $\Omega = A_{r_1,r_2}$ for $r_1,r_2\in\mathbb R$ satisfying $0<r_1<r_2$. Consider the problem \eqref{eq:pb} in $n$-dimensional spherical coordinates $(r,\theta_1,\cdots,\theta_{n-1})$, where $r\in [0,+\infty)$, $\theta_i\in[0,\pi]$ and $\theta_{n-1}\in[0,2\pi)$. Let $u(r,\theta_1,\cdots,\theta_{n-1}) = \phi(x)$ and $\tilde{m}(r,\theta_1,\cdots,\theta_{n-1})=m(x)$. 

Radial symmetry implies $u(r,\theta_1,\cdots,\theta_{n-1}) =u(r)$ and $\tilde{m}(r,\theta_1,\cdots,\theta_{n-1})=\tilde{m}(r)$. 

By change of variables, $r\in \Omega_r:=(r_1,r_2)\subseteq \mathbb R$, and the problem \eqref{eq:pb} becomes:
\begin{equation}\label{eq:pb-v3}
	\left\{
	\begin{aligned}
		u''(r) + \frac{n-1}{r} u'(r)  &+ \lambda \tilde{m}(r) u(r) =0 \ \text{in} \ \Omega_r\\
		u'(r_1) &=\beta u(r_1) \\
		u'(r_2) &= -\beta u(r_2)
	\end{aligned}
	\right.
\end{equation}

We would like to use another change of variable $r=r(t)$ to rewrite the problem \eqref{eq:pb-v3} into the preferred 1-dimensional form addressed in Lemma \ref{lem:1d}. For that purpose, we use the variational formulation of \eqref{eq:pb-v3} and get the expression of eigenvalue by Rayleigh quotient:
\begin{equation}\label{eq:spherical-rayleigh}
	\lambda = \inf_{u\in \mathcal{U}} \frac{\int_{\Omega_r} (u')^2 r^{n-1}dr + \beta\int_{\partial\Omega} u^2 r^{n-1} d\sigma}{\int_{\Omega_r} mu^2r^{n-1}dr},
\end{equation}
where $\mathcal{U}:=\left\{u\in H^1(\Omega); \int_{\Omega_r} mu^2r^{n-1}dr>0\right\}$. It is apparent that the major difference between \eqref{eq:spherical-rayleigh} and Theorem \ref{thm:principal} lies in the first term in the numerator. If we let $r=r(t)$ be a monotone increasing function on some domain $\Omega_t\subseteq \mathbb R$, then
$$\int_{r_1}^{r_2} (u')^2r^{n-1}dr = \int_{r^{-1}(r_1)}^{r^{-1}(r_2)} [u'(t)]^2 \frac{r^{n-1}(t)}{r'(t)}dt.$$
It is therefore natural to consider the solution to the ordinary differential equation
\begin{equation}
	r'(t)-r^{n-1}(t) = 0,
\end{equation}
which gives $r(t)=\gamma e^t$ for some constant $\gamma>0$ when $n=2$ and $r(t) = \left[(2-n)(t+\gamma')\right]^{\frac{1}{2-n}}$ for some constant $\gamma'\in \mathbb R$. Without loss of generality, let $r(t) = e^t$ for the planar case, and $r(t)=\left[(2-n)t\right]^{\frac{1}{2-n}}$ for higher-dimensional cases. 

\subsection{2-dimensional case}

In this section we complete the proof of Theorem \ref{thm:shell-2d}. 

\begin{proof}[Proof of Theorem \ref{thm:shell-2d}]

When $n=2$, the problem \eqref{eq:pb-v3} becomes:

\begin{equation}\label{eq:pb-v3-2d}
	\left\{
	\begin{aligned}
		u''(r) + \frac{1}{r} u'(r)  &+ \lambda \tilde{m}(r) u(r) =0 \ \text{in} \ \Omega_r\\
		u'(r_1) &=\beta u(r_1) \\
		u'(r_2) &= -\beta u(r_2)
	\end{aligned}
	\right.
\end{equation}

Let $r=e^t$ for $t\in \Omega_t:=(\ln(r_1),\ln(r_2))$, $v(t):=u(r)$, $\bar{m}(t) = \tilde{m}(r)$, and the problem \eqref{eq:pb-v3-2d} becomes

\begin{equation}\label{eq:pb-v3-2d-t}
	\left\{
	\begin{aligned}
		v''(t) + \lambda e^{2t}\bar{m}(t)v(t)=0\ \text{in}\ \Omega_t\\
		v'(\ln(r_1)) = \beta r_1 v(\ln(r_1))\\
		v'(\ln(r_2)) = -\beta r_2 v(\ln(r_2))
	\end{aligned}
	\right.
\end{equation}

Now the new weight function $M(t):=e^{2t}\bar{m}(t)$ needs scaling to fit in some $\mathcal M_{m_0,\kappa}$. 

Notice that $-r_2^2\le M(t) \le r_2^2\kappa$ and

\begin{align*}
	\int_{\ln(r_1)}^{\ln(r_2)} M(t) dt =\int_{r_1}^{r_2} r\tilde{m}(r)dr\\
	= \frac{1}{2\pi} \int_{\Omega} m(x)dx \le -\frac{m_0}{2}(r_2^2-r_1^2),
\end{align*}

so we can take a real constant $q$ satisfying 
$$q> \frac{m_0}{2(\ln(r_2)-\ln(r_1))},$$
so that 
$$\int_{\ln(r_1)}^{\ln(r_2)} M(t) dt \le -\frac{m_0(r_2^2-r_1^2)}{2q(\ln(r_2)-\ln(r_1))}|\Omega_t|.$$

Then let $\bar{M}(t):=M(t)/r_2^2$, and clearly $-1\le \bar{M}(t)\le \kappa$ and
$$\int_{\ln(r_1)}^{\ln(r_2)} \bar{M}(t)dt \le - m_0'|\Omega_t|,$$
where 
$$m_0' := \frac{m_0(r_2^2-r_1^2)}{2qr_2^2(\ln(r_2)-\ln(r_1))} \in (0,1).$$

So $\bar{M}\in \mathcal{M}_{m_0',\kappa}$. Letting $\beta' := \beta r_1$ and $\lambda' := r_2^2\lambda$, the problem \eqref{eq:pb-v3-2d-t} becomes:

\begin{equation}\label{eq:pb-v3-2d-t2}
	\left\{
	\begin{aligned}
		v''(t) + \lambda' \bar{M}(t)v(t)=0\ \text{in}\ \Omega_t\\
		v'(\ln(r_1)) = \beta' v(\ln(r_1))\\
		v'(\ln(r_2)) = -\beta' r_1^{-1} r_2 v(\ln(r_2))
	\end{aligned}
	\right.
\end{equation}

Since the problem \ref{eq:pb-v3-2d-t2} differs from the original problem \eqref{eq:pb} only at a boundary point, we can still apply Lemma \ref{lem:1d} to it, as the result will only differ by a set of measure zero. So the optimal weight $\bar{M}^*$ for \eqref{eq:pb-v3-2d-t2} satisfies
$$\bar{M}^* = \kappa \chi_{E_t^*} - \chi_{\Omega_t\setminus E_t^*}\ \text{a.e.}\ t\in \Omega_t$$
for some optimal set $E_t^*$ given by Lemma \ref{lem:1d}. More specifically, let $c':=\frac{1-m_0'}{1+\kappa}$ and $\beta^*=\beta^*(c',\kappa)$ as defined in Theorem \ref{thm:sol1d}, then we have:
\begin{itemize}
	\item if $\beta' > \frac{\beta^*}{\ln(r_2)-\ln(r_1)}$, we have
	$$E_t^* = \left(\frac{(1+c')\ln(r_1) + (1-c')\ln(r_2)}{2}, \frac{(1-c')\ln(r_1)+(1+c')\ln(r_2)}{2}\right);$$
	\item if $\beta' = \frac{\beta^*}{\ln(r_2)-\ln(r_1)}$, $E_t^*$ can be any open interval of length $c'(\ln(r_2)-\ln(r_1))$; 
	\item if $\beta' < \frac{\beta^*}{\ln(r_2)-\ln(r_1)}$, then $E_t^*=(\ln(r_1), (1-c')\ln(r_1)+c'\ln(r_2))$ or $E_t^*=(c'\ln(r_1)+(1-c')\ln(r_2), \ln(r_2))$. 
\end{itemize} 

Using $r=e^t$, we get the optimal set $E_r^*$ for the variable $r$. The proof of Theorem \ref{thm:shell-2d} is thus completed.  

\end{proof}

\subsection{$n$-dimensional case for $n\ge 3$}

In this section, we complete the proof of Theorem \ref{thm:shell-3d}. 

\begin{proof}[Proof of Theorem \ref{thm:shell-3d}]

In such case, we let $r=[(2-n)t]^{\frac{1}{2-n}}$ for 
$$t\in \Omega_t:=\left(\frac{r_1^{2-n}}{2-n}, \frac{r_2^{2-n}}{2-n}\right)\subseteq (-\infty, 0).$$

Above all, notice that
$$r'(t)=\left[(2-n)t\right]^{\frac{n-1}{2-n}}>0$$
whenever $t<0$, hence $r(t)$ is monotone increasing in $\Omega_t\subseteq \mathbb R$. 

Now let $v(t):=u(r)$, $\bar{m}(t)=\tilde{m}(r)$, and the problem \eqref{eq:pb-v3} becomes

\begin{equation}\label{eq:pb-v3-t}
	\left\{
	\begin{aligned}
		v''(t) + \lambda [(2-n)t]^{\frac{2n-2}{2-n}} \bar{m}(t)v(t)=0\ \text{in}\ \Omega_t\\
		v'\left(\frac{r_1^{2-n}}{2-n}\right) = \beta r_1^{n-1} v\left(\frac{r_1^{2-n}}{2-n}\right)\\
		v'\left(\frac{r_2^{2-n}}{2-n}\right)=- \beta r_2^{n-1} v\left(\frac{r_2^{2-n}}{2-n}\right)
	\end{aligned}
	\right.
\end{equation}

Now consider the new weight function $M(t) = g(t)\bar{m}(t)$, where
$$g(t)=[(2-n)t]^{\frac{2n-2}{2-n}}.$$
Then 
$$g'(t) = (2n-2)[(2-n)t]^{\frac{3n-4}{2-n}}.$$
Notice that, when $t\in \Omega_t$, $t<0$, and we have $g'(t)>0$, so $g(t)$ is monotone increasing in $\Omega_t$. Therefore, 
$$r_1^{2n-2}=g\left(\frac{r_1^{2-n}}{2-n}\right)\le g(t)\le g\left(\frac{r_2^{2-n}}{2-n}\right)=r_2^{2n-2},$$
so $-r_2^{2n-2}\le M(t)\le r_2^{2n-2}\kappa$. Moreover, 

$$\int_{\frac{r_1^{2-n}}{2-n}}^{\frac{r_2^{2-n}}{2-n}} M(t) dt = \int_{r_1}^{r_2} r^{n-1} \tilde{m}(r) dr.$$

On the other hands, 

\begin{align*}
	\int_{\Omega} m(x) dx  &= \int_0^{2\pi} \int_{0}^{\pi}\cdots \int_0^{\pi} \int_{r_1}^{r_2} \tilde{m}(r) r^{n-1} \prod_{k=1}^{n-2} \sin^{n-k-1}(\theta_k) dr d\theta_1\cdots d\theta_{n-2}d\theta_{n-1}\\
	&=2\pi \prod_{k=1}^{n-2} \int_{0}^{\pi} \sin^k(x)dx \int_{r_1}^{r_2} \tilde{m}(r)r^{n-1}dr.
\end{align*}

Using the Gamma function, we have
\begin{equation}
	\int_\Omega m(x)dx = 2\pi ^{\frac{n}{2}} \prod_{k=1}^{n-2} \frac{\Gamma\left(\frac{k+1}{2}\right)}{\Gamma\left(\frac{k}{2} + 1\right)} \int_{r_1}^{r_2} \tilde{m}(r)r^{n-1}dr.
\end{equation}

Therefore, 
\begin{align*}
	\int_{\Omega_t} M(t)dt &= \left( 2\pi ^{\frac{n}{2}} \prod_{k=1}^{n-2} \frac{\Gamma\left(\frac{k+1}{2}\right)}{\Gamma\left(\frac{k}{2} + 1\right)}\right)^{-1} \int_\Omega m(x) dx \\
	&\le-m_0|\Omega|\left( 2\pi ^{\frac{n}{2}} \prod_{k=1}^{n-2} \frac{\Gamma\left(\frac{k+1}{2}\right)}{\Gamma\left(\frac{k}{2} + 1\right)}\right)^{-1}.
\end{align*}

Since 
$$|\Omega| = |B_{r_2}(0)|-|B_{r_1}(0)|=\frac{\pi^{\frac{n}{2}}\left(r_2^2-r_1^2\right)}{\Gamma\left(\frac{n}{2} + 1\right)},$$
we have
\begin{equation}
	\int_{\Omega_t} M(t) dt \le -\frac{m_0\left(r_2^2-r_1^2\right)}{2\Gamma\left(\frac{n}{2}+1\right)} \prod_{k=1}^{n-2}\frac{\Gamma\left(\frac{k}{2}+1\right)}{\Gamma\left(\frac{k+1}{2}\right)}.
\end{equation}

Now it suffices to choose a positive constant $q$ satisfying
$$q > \max\left\{m_0, \frac{2|\Omega_t|r_2^{2n-2}\Gamma\left(\frac{n}{2}+1\right)}{r_2^2-r_1^2} \prod_{k=1}^{n-2}\frac{\Gamma\left(\frac{k+1}{2}\right)}{\Gamma\left(\frac{k}{2}+1\right)}\right\}$$
and let $\bar{M}(t):=r_2^{2-2n}M(t)$ and $m_0':=q^{-1}m_0$, then $\bar{M}(t)\in[-1,\kappa]$ and
$$\int_{\Omega_t} \bar{M}(t)dt \le -m_0' |\Omega_t|.$$

So $\bar{M}\in \mathcal{M}_{m_0',\kappa}$. Letting $\beta'=\beta r_1^{n-1}$ and $\lambda' = r_2^{2n-2}\lambda$, the problem \eqref{eq:pb-v3-t} becomes:

\begin{equation}\label{eq:pb-v3-t2}
	\left\{
	\begin{aligned}
		v''(t) + \lambda' \bar{M}(t)v(t) &=0\ \text{in}\ \Omega_t\\
		v'\left(\frac{r_1^{2-n}}{2-n}\right)&= \beta'  v\left(\frac{r_1^{2-n}}{2-n}\right)\\
		v'\left(\frac{r_2^{2-n}}{2-n}\right) &= -\beta' \left(\frac{r_2}{r_1}\right)^{n-1} v\left(\frac{r_2^{2-n}}{2-n}\right)
	\end{aligned}
	\right.
\end{equation}

Since the problem \eqref{eq:pb-v3-t2} differs from the original problem \eqref{eq:pb} only at a boundary point, we can still apply Lemma \ref{lem:1d} to it, as the result will only differ by a set of measure zero. So the optimal weight $\bar{M}^*$ for \eqref{eq:pb-v3-t2} satisfies
$$\bar{M}^* = \kappa \chi_{E_t^*} - \chi_{\Omega_t\setminus E_t^*}\ \text{a.e.}\ t\in \Omega_t$$
for some optimal set $E_t^*$ given by Lemma \ref{lem:1d}. More specifically, let $c':=\frac{1-m_0'}{1+\kappa}$ and $\beta^*=\beta^*(c',\kappa)$ as defined in Theorem \ref{thm:sol1d}, then we have:
\begin{itemize}
	\item if $\beta' > \frac{\beta^*(2-n)}{r_2^{2-n}-r_1^{2-n}}$, 
	$$E_t^*=\left(\frac{(1+c')r_1^{2-n}+(1-c')r_2^{2-n}}{4-2n}, \frac{(1-c')r_1^{2-n}+(1+c')r_2^{2-n}}{4-2n}\right);$$
	\item if $\beta' = \frac{\beta^*(2-n)}{r_2^{2-n}-r_1^{2-n}}$, 
	$E_t^*$ can be any open interval of length $\frac{c'(r_2^{2-n}-r_1^{2-n})}{2-n}$; 
	\item if $\beta' < \frac{\beta^*(2-n)}{r_2^{2-n}-r_1^{2-n}}$, we have 
	$$E_t^*=\left(\frac{r_1^{2-n}}{2-n}, \frac{(1-c')r_1^{2-n} + c'r_2^{2-n}}{2-n}\right)$$ or 
	$$E_t^*=\left(\frac{c'r_1^{2-n}+(1-c')r_2^{2-n}}{2-n}, \frac{r_2^{2-n}}{2-n}\right).$$
\end{itemize}

Using $r=[(2-n)t]^{\frac{1}{2-n}}$, we get the optimal set $E_r^*$ for the variable $r$. The proof of Theorem \ref{thm:shell-3d} is thus completed.  

\end{proof}

%% file: sections/acknowledgment.tex
\section*{Acknowledgments}\label{sec:acknowledgement}
This work is partially supported by the project SGS01/PřF/2024 ``Aspects of algebra, mathematical analysis and number theory". 

%% file: bibliography/references.bib
@Article{Afrouzi1999,
  author  = {Afrouzi, G. A. and Brown, K. J.},
  journal = {Proceedings of the American Mathematical Society},
  title   = {On principal eigenvalues for boundary value problems with indefinite weight and {R}obin boundary conditions},
  year    = {1999},
  issn    = {1088-6826},
  number  = {1},
  pages   = {125--130},
  volume  = {127},
  doi     = {10.1090/s0002-9939-99-04561-x},
}

@Article{Bocher1914,
  author  = {Bôcher, Maxime},
  journal = {Bulletin of the American Mathematical Society},
  title   = {The smallest characteristic numbers in a certain exceptional case},
  year    = {1914},
  issn    = {1088-9485},
  number  = {1},
  pages   = {6--9},
  volume  = {21},
  doi     = {10.1090/s0002-9904-1914-02560-1},
}

@Article{Lou2006,
  author  = {Lou, Yuan and Yanagida, Eiji},
  journal = {Japan Journal of Industrial and Applied Mathematics},
  title   = {Minimization of the principal eigenvalue for an elliptic boundary value problem with indefinite weight, and applications to population dynamics},
  year    = {2006},
  issn    = {1868-937X},
  number  = {3},
  pages   = {275--292},
  volume  = {23},
  doi     = {10.1007/bf03167595},
}

@Article{Lamboley2016,
  author  = {Lamboley, Jimmy and Laurain, Antoine and Nadin, Grégoire and Privat, Yannick},
  journal = {Calculus of Variations and Partial Differential Equations},
  title   = {Properties of optimizers of the principal eigenvalue with indefinite weight and {R}obin conditions},
  year    = {2016},
  issn    = {1432-0835},
  number  = {6},
  volume  = {55},
  doi     = {10.1007/s00526-016-1084-6},
}

@Article{Jiang2023,
  author  = {Mengping Jiang},
  journal = {Advances in Applied Mathematics},
  title   = {Research on the minimization problem of principal eigenvalue in spherical shell domain},
  year    = {2023},
  issn    = {2324-8009},
  number  = {09},
  pages   = {3826--3833},
  volume  = {12},
  doi     = {10.12677/aam.2023.129376},
}

@Article{Hintermueller2011,
  author  = {Hintermüller, M. and Kao, C.-Y. and Laurain, A.},
  journal = {Applied Mathematics and Optimization},
  title   = {Principal eigenvalue minimization for an elliptic problem with indefinite weight and {R}obin boundary conditions},
  year    = {2011},
  issn    = {1432-0606},
  number  = {1},
  pages   = {111--146},
  volume  = {65},
  doi     = {10.1007/s00245-011-9153-x},
}

@Article{Kao2008,
  author  = {Chiu-Yen Kao and Yuan Lou and Eiji Yanagida},
  journal = {Mathematical Biosciences and Engineering},
  title   = {Principal eigenvalue for an elliptic problem with indefinite weight on cylindrical domains},
  year    = {2008},
  issn    = {1551-0018},
  number  = {2},
  pages   = {315--335},
  volume  = {5},
  doi     = {10.3934/mbe.2008.5.315},
}

@Article{Brown1980,
  author  = {K. J. Brown and S. S. Lin},
  journal = {Journal of Mathematical Analysis and Applications},
  title   = {On the existence of positive eigenfunctions for an eigenvalue problem with indefinite weight function},
  year    = {1980},
  issn    = {0022-247X},
  number  = {1},
  pages   = {112--120},
  volume  = {75},
  doi     = {10.1016/0022-247x(80)90309-1},
}

@Article{Senn1982,
  author  = {Senn, Stefan and Hess, Peter},
  journal = {Mathematische Annalen},
  title   = {On positive solutions of a linear elliptic eigenvalue problem with neumann boundary conditions},
  year    = {1982},
  issn    = {1432-1807},
  number  = {4},
  pages   = {459--470},
  volume  = {258},
  doi     = {10.1007/bf01453979},
}

@Article{MunozCastaneda2015,
  author  = {Muñoz-Castañeda, J. M. and Mateos Guilarte, J.},
  journal = {Physical Review D},
  title   = {$\delta$-$\delta'$ generalized {R}obin boundary conditions and quantum vacuum fluctuations},
  year    = {2015},
  issn    = {1550-2368},
  number  = {2},
  pages   = {025028},
  volume  = {91},
  doi     = {10.1103/physrevd.91.025028},
}

@Article{Robin2021,
  author  = {Robin, Caroline and Savage, Martin J. and Pillet, Nathalie},
  journal = {Physical Review C},
  title   = {Entanglement rearrangement in self-consistent nuclear structure calculations},
  year    = {2021},
  issn    = {2469-9993},
  number  = {3},
  pages   = {034325},
  volume  = {103},
  doi     = {10.1103/physrevc.103.034325},
}

@Article{Fritze2005,
  author  = {Fritze, Denny and Marburg, Steffen and Hardtke, Hans-Jürgen},
  journal = {Computers and Structures},
  title   = {{FEM}–{BEM}-coupling and structural–acoustic sensitivity analysis for shell geometries},
  year    = {2005},
  issn    = {0045-7949},
  number  = {2–3},
  pages   = {143--154},
  volume  = {83},
  doi     = {10.1016/j.compstruc.2004.05.019},
}

@Article{Clarte2021,
  author  = {Clarté, Thibaut T. and Schaeffer, Nathanaël and Labrosse, Stéphane and Vidal, Jérémie},
  journal = {Journal of Fluid Mechanics},
  title   = {The effects of a {R}obin boundary condition on thermal convection in a rotating spherical shell},
  year    = {2021},
  issn    = {1469-7645},
  volume  = {918},
  doi     = {10.1017/jfm.2021.356},
}

@Article{Esposito1999,
  author  = {Esposito, Giampiero and Yu. Kamenshchik, Alexander and Kirsten, Klaus},
  journal = {International Journal of Modern Physics A},
  title   = {Zero-point energy of a conducting spherical shell},
  year    = {1999},
  issn    = {1793-656X},
  number  = {02},
  pages   = {281--300},
  volume  = {14},
  doi     = {10.1142/s0217751x99000154},
}
